\documentclass[10pt]{amsart}
\usepackage{amssymb,amsfonts}

\newcommand{\CC}{{\mathbb C}}
\newcommand{\RR}{{\mathbb R}}
\newcommand{\ZZ}{{\mathbb Z}}
\newcommand{\n}{\nabla}
\renewcommand{\phi}{\varphi}
\renewcommand{\epsilon}{\varepsilon}
\newcommand{\half}{\tfrac{1}{2}}
\renewcommand{\div}{\operatorname{div}}
\newcommand{\grad}{\operatorname{grad}}

\newcommand{\Spin}{\operatorname{Spin}}
\newcommand{\Pin}{\operatorname{Pin}}

\newcommand{\scal}{\operatorname{scal}}
\newcommand{\tr}{\operatorname{tr}}

\newcommand{\Rinv}{{R^{\rm{inv}}} } 
\newcommand{\Rpsc}{{R^{\rm{psc}}} } 
\newcommand{\ind}{{\operatorname{ind}}}
\newcommand{\vol}{{\operatorname{vol}}}

\newtheorem{thm}{Theorem}[section]
\newtheorem{lemma}[thm]{Lemma}
\newtheorem{prop}[thm]{Proposition}
\newtheorem{cor}[thm]{Corollary}
\newtheorem{remark}[thm]{Remark}
\newtheorem{remarks}[thm]{Remarks}
\newtheorem{definition}[thm]{Definition}
\newtheorem{notation}[thm]{Notation}
\newtheorem{example}[thm]{Example}

\begin{document}

\title
[On the space of metrics with invertible Dirac operator]
{On the space of metrics with invertible Dirac operator}

\author{Mattias Dahl}

\address{
Institutionen f\"or Matematik\\
Kungl Tekniska H\"ogskolan\\
100 44 Stockholm\\
Sweden
}
 
\email{dahl@math.kth.se}



\date{\today}

\begin{abstract}
On a compact spin manifold we study the space of Riemannian metrics
for which the Dirac operator is invertible. The first main result is 
a surgery theorem stating that such a metric can be extended over the
trace of a surgery of codimension at least three. We then prove that
if non-empty the space of metrics with invertible Dirac operators is
disconnected in dimensions $n \equiv  0,1,3,7 \mod 8$, $n \geq
5$. As a corollary follows results on the existence of metrics with
harmonic spinors by Hitchin and B\"ar. Finally we use computations of 
the eta invariant by Botvinnik and Gilkey to find metrics with
harmonic spinors on simply connected manifolds with a cyclic group
action. In particular this applies to spheres of all dimensions $n
\geq 5$.
\end{abstract}

\maketitle

\section{Introduction}

Let $(M,g)$ be a Riemannian spin manifold, we will always assume that
such a manifold comes equipped with an of orientation and a 
spin structure. We denote by $M^-$ the same manifold with the opposite 
orientation. The Dirac operator $D^g$ is a first order elliptic
differential operator acting on smooth sections of the spinor bundle
$\Sigma M$. If $M$ has a boundary we will only consider Riemannian
metrics on $M$ which have a product structure in a neighbourhood of
the boundary.

For a Riemannian manifold $(M,g)$ with boundary $\partial M$ we denote 
by $(M_{\infty},g)$ the same manifold with the half-infinite cylinder 
$([0,\infty) \times \partial M, dx^2 + g|_{ \partial M})$ 
attached along the boundary (here we abuse notation slightly by using 
the same symbol $g$ for the metric on $M$ and the metric on $M_{\infty}$).
If $M$ is closed, that is compact with no boundary, we set 
$(M_{\infty},g) = (M,g)$.

We denote by $C_0^{\infty}(\Sigma M)$ the space of compactly supported
smooth sections of $\Sigma M$. On a complete Riemannian manifold
$(M,g)$ we denote by $L^2(\Sigma M)$ and $H^1(\Sigma M)$ the
completions of $C_0^{\infty}(\Sigma M)$ with respect to the $L^2$-norm
$\| \cdot \|$ and the first Sobolev norm $\| \cdot \|_1$. 

If $(M,g)$ is compact without boundary the operator $D^g$ has a
self-adjoint extension to $L^2(\Sigma M)$ with domain $H^1(\Sigma M)$.
This is a Fredholm operator with discrete spectrum
\cite[Chap. 3, \S 5]{Lawson_Michelsohn89}. If $(M,g)$ is compact with
non-empty boundary we consider the Dirac operator $D^g$ on the
manifold $(M_{\infty},g)$ with cylindrical ends. In this case we also
have a self-adjoint extension to $L^2(\Sigma M_{\infty})$ with domain
$H^1(\Sigma M_{\infty})$, see \cite[Sec. 3.6.2]{Bleecker_Booss-Bavnbek03}.

Now suppose $(M,g)$ is compact, possibly with boundary. The operator 
$D^g$ is invertible with a bounded inverse if and only if it has a 
spectral gap around $0$, that is if there is an $\epsilon>0$ such that 
$\|D^g \phi \|^2 \geq \epsilon \| \phi \|^2$ for all 
$\phi \in L^2(\Sigma M_{\infty})$.

\begin{definition} \label{def_of_Rinv}
Suppose $M$ is a compact spin manifold. We define $\Rinv(M)$ to be 
the set of Riemannian metrics $g$ on $M$ for which $D^g$ is invertible 
with a bounded inverse. 
\end{definition}

Let $R(M)$ be set of all Riemannian metrics on $M$. If $M$ is a closed 
spin manifold then $\Rinv(M)$ is an open subset of $R(M)$ in the 
$C^1$-topology and if non-empty $\Rinv(M)$ is dense in $R(M)$ in the 
$C^k$-topology for all $k \geq 1$, see \cite[Prop. 3.2]{Baer_Dahl02}.

\begin{prop} \label{prop:inv_restriction_to_bdry}
If $g \in \Rinv(M)$ then $g|_{ \partial M} \in \Rinv(\partial M)$.
\end{prop}

\begin{proof}
Suppose that the Dirac operator for $g|_{ \partial M}$ is not invertible. 
Then there is a $\phi \neq 0$ such that $D^{g|_{ \partial M}} \phi = 0$. 
If we extend $\phi$ to the cylindrical end of $(M_{\infty},g)$ by parallel 
transport in the normal direction and then multiply with a cut-off 
function having small gradient we can construct compactly supported 
$\psi$ on $M_{\infty}$ for which $\|\psi\|^2 = 1$ and $\|D^g \psi\|^2$ is 
arbitrarily small. 
\end{proof}

\begin{definition}
Let $M, N$ be compact spin manifolds without boundary. 
\begin{enumerate}
\item Metrics $g^0, g^1 \in \Rinv(M)$ are called isotopic
if there is a smooth path of metrics $g_t \in \Rinv(M)$, $t \in \RR$,
such that $g_t =g^0$ for $t \leq 0$ and  $g_t =g^1$ for $t \geq 1$. 
\item Metrics $g^0, g^1 \in \Rinv(M)$ are called concordant if there is a
metric $\overline{g} \in \Rinv([0,1] \times M)$ such that 
$\overline{g}|_{ \{ i \} \times M} = g^i$, $i=0,1$. 
\item Metrics $g^0 \in \Rinv(M)$, $g^1 \in \Rinv(N)$, are called 
bordant if there is a manifold $W$ and a metric $g^W \in \Rinv(W)$ 
so that $\partial(W,g^W) = (M,g^0) \sqcup (N^-,g^1)$.
\end{enumerate}
\end{definition}
The Dirac operator is intimately related to the scalar curvature.
\begin{definition}
Let $\Rpsc(M)$ be the set of Riemannian metrics on $M$ with positive 
scalar curvature. 
\end{definition}

From the Lichnerowicz formula 
$(D^g)^2 = (\n^g)^* \n^g + \frac{1}{4} \scal^g$ it follows that 
$\Rpsc(M) \subset \Rinv(M)$. There are the corresponding relations 
psc-isotopic/psc-concordant/psc-bordant for metrics in $\Rpsc(M)$. The 
Lichnerowicz formula implies that if two metrics are 
psc-isotopic/psc-concordant/psc-bordant then they are 
isotopic/concordant/bordant.

The main idea of this paper is to study the space $\Rinv(M)$ using 
techniques from the study of $\Rpsc(M)$. In Section 
\ref{Sec:Constructions} we will look at ways of constructing
Riemannian manifolds with invertible Dirac operator, the most powerful
of which will be the extension of a metric with invertible Dirac
operator to the trace of a surgery of codimension at least three. The
main result of the paper is in Section \ref{Sec:Detecting} where we
use the Index Theorem to detect non-concordant metrics in $\Rinv(M)$
in dimensions $n \equiv  0,1,3,7 \mod 8$, $n \geq 5$. The construction
of these non-concordant metrics uses known examples of ``exotic''
metrics in $\Rpsc(S^n)$ which do not bound metrics in
$\Rpsc(D^{n+1})$. This result shows that if non-empty $\Rinv(M)$ is
disconnected, which unifies and strengthens results by Hitchin
and B\"ar on the existence of metrics with non-trivial harmonic
spinors. Finally in Section \ref{Sec:Cyclic} we leave the study of
$\Rinv(M)$. Instead we use computations of the eta invariant by
Botvinnik and Gilkey to find metrics with harmonic spinors on simply
connected manifolds with a cyclic group action. In particular we find
metrics with harmonic spinors on spheres of all dimensions $n \geq 5$.  

\section{Constructions} 
\label{Sec:Constructions}

In this section we will study three constructions of new Riemannian 
manifolds with invertible Dirac operators from old ones.

\subsection{Attaching isometric boundary components}

Let $M$ be a manifold with boundary $\partial M$. Suppose that the 
boundary is a disjoint union 
$\partial M = \partial^+ M \sqcup \partial^- M \sqcup \partial^0 M$
where $\partial^+ M \cong N$ and $\partial^- M \cong N^-$ for some compact
spin manifold $N$ and where $\partial^0 M$ might be empty. 

Suppose $g \in \Rinv(M)$ is such that 
$g|_{\partial^+ M} = g|_{\partial^- M} = h$ for some metric $h$ on $N$.
For $t > 0$ let $(M',g'_t)$ be $(M,g)$ with the 
cylinder $([0,t] \times N, dx^2 + h)$ attached by $\{ 0 \} \times N$
along $\partial^+ M$ and by $\{ t \} \times N$ along $\partial^- M$. 
The manifolds $M'$ depend on $t$ but are all diffeomorphic so we
identify them.

\begin{prop} \label{prop:attachingboundaries}
Let $(M',g'_t)$ be constructed from $(M,g)$ as above.
Then there is $T>0$ so that $g'_t \in \Rinv(M')$ for 
all $t>T$.
\end{prop}

Note that the manifold $M$ is not assumed to be connected. 

\begin{proof}
Since $g \in \Rinv(M)$ there is $\epsilon^g > 0$ so that 
$\| D^g \phi \|^2 \geq \epsilon^g \| \phi \|^2$ for all 
$\phi\in L^2(\Sigma M_{\infty})$. Set 
$\epsilon = \frac{\epsilon^g}{8}$ and choose $T>0$ so that $t>T$ 
implies $\frac{6}{t^2} \leq  \frac{\epsilon^g}{8}$.

Let $t>T$ and take $\phi \in C_0^{\infty} (\Sigma M'_{\infty})$. Let 
$\chi: [0,t] \times N \to [0,1]$ be a smooth function such that 
$\chi = 1$ near $\{ 0 \} \times N$, $\chi = 0$ near $\{ t \} \times N$, 
and $|\grad \chi| \leq \frac{2}{t}$. A straight-forward computation 
shows that 
$$
|D\phi|^2 \geq
\tfrac{1}{2} |D (\chi \phi)|^2 + \tfrac{1}{2} |D ((1-\chi) \phi)|^2
-\tfrac{3}{2} |\grad \chi|^2 |\phi|^2 
$$
so
$$
\| D\phi \|^2_{[0,t] \times N} \geq
\tfrac{1}{2} \| D ( \chi \phi) \|^2_{[0,t] \times N}
+
\tfrac{1}{2} \| D ( (1-\chi) \phi) \|^2_{[0,t] \times N}
-
\frac{6}{t^2}  \| \phi \|^2_{[0,t] \times N}.
$$
We define the spinor field $\psi \in C_0^{\infty} (M_{\infty})$ as 
follows. On $M$ and on $[0,\infty) \times \partial^0 M$ we set 
$\psi = \phi$. At $\partial^+ M$ we first attach 
$([0,t] \times N, dx^2 + h)$ along $\{ 0 \} \times N$ and set 
$\psi = \chi \phi$ on this piece, followed by $\psi=0$ on the 
half-infinite cylinder attached along $\{ t \} \times N$. In the 
same way we attach $([0,t] \times N, dx^2 + h)$ at $\partial^- M$ 
along $\{ t \} \times N$ and there we set $\psi = (1-\chi) \phi$ 
followed by $\psi=0$ on the half-infinite cylinder attached along 
$\{ 0 \} \times N$. Using the above estimate we get
\begin{equation*}
\begin{split}
\| D\phi \|^2_{M'_{\infty}}
&=
\| D\phi \|^2_{M} + \| D\phi \|^2_{[0,\infty) \times \partial^0 M}
+ \| D\phi \|^2_{[0,t] \times N} \\
&\geq
\| D\phi \|^2_{M} + \| D\phi \|^2_{[0,\infty) \times \partial^0 M} \\
&\qquad{}+ 
\tfrac{1}{2} \| D ( \chi \phi) \|^2_{[0,t] \times N}
+
\tfrac{1}{2} \| D ( (1-\chi) \phi) \|^2_{[0,t] \times N}
-
\frac{6}{t^2}  \| \phi \|^2_{[0,t] \times N} \\
&\geq
\tfrac{1}{2} \left(
\| D\phi \|^2_{M} + \| D\phi \|^2_{[0,\infty) \times \partial^0 M}
\right. \\
&\qquad \left. {}
+ \| D ( \chi \phi) \|^2_{[0,t] \times N}
+ \| D ( (1-\chi) \phi) \|^2_{[0,t] \times N} \right)
- \frac{6}{t^2}  \| \phi \|^2_{M'_{\infty}} \\
&=
\tfrac{1}{2} \| D\psi \|^2_{M_{\infty}} 
- \frac{6}{t^2}  \| \phi \|^2_{M'_{\infty}} \\
&\geq
\frac{\epsilon^g}{2} \| \psi \|^2_{M_{\infty}} 
- \frac{6}{t^2}  \| \phi \|^2_{M'_{\infty}} \\
&=
\frac{\epsilon^g}{2} \left(
\| \phi \|^2_{M} + \| \phi \|^2_{[0,\infty) \times \partial^0 M}
+ \| \chi \phi \|^2_{[0,t] \times N}
+ \| (1-\chi) \phi \|^2_{[0,t] \times N} \right) \\
&\qquad{} -
\frac{6}{t^2}  \| \phi \|^2_{M'_{\infty}} \\
&\geq
\frac{\epsilon^g}{2} \left(
\| \phi \|^2_{M} + \| \phi \|^2_{[0,\infty) \times \partial^0 M}
+ \tfrac{1}{2} \| \phi \|^2_{[0,t] \times N} \right)
- \frac{6}{t^2}  \| \phi \|^2_{M'_{\infty}} \\
&\geq
\left(\frac{\epsilon^g}{4} - \frac{6}{t^2} \right)
\| \phi \|^2_{M'_{\infty}} \\
&\geq
\epsilon \| \phi \|^2_{M'_{\infty}}.
\end{split}
\end{equation*}
Since $C_0^{\infty}(\Sigma M'_{\infty})$ is dense in 
$L^2(\Sigma M'_{\infty})$ this shows that $g'_t \in \Rinv(M')$.
\end{proof}

Proposition \ref{prop:attachingboundaries} has the following 
corollary.

\begin{cor} \label{cor:conc&bord_eqrel}
Concordance and bordance are equivalence relations.
\end{cor}

\subsection{Generalized cylinders}

Let $M$ be a compact spin manifold of dimension $n$ and let $g_{\tau}$ 
be a smooth curve of metrics on $M$ parametrized by $\tau \in I$, where
$I$ is an interval. The product $\overline{M} = I \times M$ equipped
with the metric $\overline{g} = d\tau^2 + g_{\tau}$ is called a
generalized cylinder over $M$. We are going to recall some facts about
the spinor bundle and the Dirac operator on a generalized
cylinder. All these facts are conveniently collected in
\cite{Baer_Gauduchon_Moroianu03}.

The spin structure on $M$ induces in a unique way a spin structure on  
$\overline{M}$. The spinor bundle on $(\overline{M},\overline{g})$ is
related to the spinor bundle on $(M,g_\tau)$ by 
$\Sigma_{(\tau,x)} \overline{M} = \Sigma_{x} M$ if $n$ is even and 
$\Sigma_{(\tau,x)}^{\pm} \overline{M} = \Sigma_{x} M$ if $n$ is odd. 
Denote Clifford multiplication on $\Sigma \overline{M}$ by $\cdot$ and 
Clifford multiplication on $\Sigma M$ by $\bullet_\tau$. If $n$ is even 
we have $X \bullet_\tau \phi = \nu \cdot X \cdot \phi$ and if $n$ is odd 
$X \bullet_\tau \phi = \pm \nu \cdot X \cdot \phi$ for 
$\phi \in \Sigma^{\pm} \overline{M}$. Here $\nu = \partial_\tau$ is the 
normal of $\{\tau\} \times M$ in $(\overline{M}, \overline{g})$.

Let $\phi$ be a section of $\Sigma \overline{M}$. The Dirac operators 
on $M$ and $\overline{M}$ are related by 
\begin{equation} \label{Diraconcylinder}
\nu \cdot D^{\overline{g}} \phi =
\left(D^{g_{\tau}}
+ \frac{n}{2} H -\n_{\nu}^{\overline{g}}
\right) \phi.
\end{equation}
Here $H$ is the mean curvature of $\{\tau\} \times M$ in 
$(\overline{M}, \overline{g})$, 
\begin{equation} \label{meancurvature}
H = - \frac{1}{2n} \tr_{g_\tau} \left(
\dot{g}_\tau \right),
\end{equation}
and if $n$ is odd the operator
$D^{g_{\tau}}$ acts on sections of $\Sigma \overline{M}$ by 
$\left(\begin{smallmatrix}
D^{g_{\tau}} & 0 \\
0 & -D^{g_{\tau}} 
\end{smallmatrix}\right)$.
Let $\dot{g}_\tau = \partial_\tau g_\tau$ and define the operator
${\mathfrak D}^{\dot{g}_{\tau}}$ by
${\mathfrak D}^{\dot{g}_{\tau}} \phi
=
\sum_{i,j=1}^n \dot{g}_{\tau} (e_i,e_j) e_i \bullet_\tau
\n_{e_j}^{g_{\tau}} \phi$
where $e_1,\dots, e_n$ is an orthonormal basis of $TM$. The commutator
of $\n_{\nu}^{\overline{g}}$ and $D^{g_{\tau}}$ is given by 
\cite[Eqn. (23)]{Baer_Gauduchon_Moroianu03}
\begin{equation} \label{commutator}
\left[ \n_{\nu}^{\overline{g}}, D^{g_{\tau}} \right]\phi
=
- \tfrac{1}{2} {\mathfrak D}^{\dot{g}_{\tau}} \phi
+ \tfrac{1}{4} \grad^{g_{\tau}} 
\left( 
\tr_{g_{\tau}} \left( \dot{g}_{\tau} \right)
\right) \bullet_\tau \phi
-
\tfrac{1}{4} \div^{g_{\tau}}
\left( \dot{g}_{t} \right) \bullet_\tau \phi.
\end{equation}

Now suppose $g_\tau$, $\tau\in [0,1]$, is a smooth curve of metrics in 
$\Rinv(M)$ with $g_\tau = g_{0}$ for $\tau$ near $0$ and 
$g_\tau = g_1$ for $\tau$ near $1$. Define metrics $\overline{g}_t$ on 
$\overline{M}_t = [0,t] \times M$ by 
$\overline{g}_t = d\tau^2 + g_{\tau/t}$ for $t > 0$. 
Since the $\overline{M}_t$ are all diffeomorphic we identify 
them as $\overline{M}$.

\begin{prop} \label{prop:stretchcylinder}
Suppose $(\overline{M},\overline{g}_t )$ is constructed from 
$M$ and $g_\tau$ as above. Then there exists $T > 0$ such that 
$\overline{g}_t \in \Rinv(\overline{M})$ for all $t > T$. 
\end{prop}

\begin{proof}
Since $g_\tau$ is defined for $\tau$ in a compact interval and since 
$g_\tau \in \Rinv(M)$ there is a constant $C > 0$ so that 
\begin{gather}
\frac{1}{C} \int_M \left|\phi \right|^2 d \vol^{g_{\tau}} 
\leq
\int_M \left| D^{g_\tau} \phi \right|^2 d \vol^{g_{\tau}},
\label{C1} \\
\left|
\tfrac{1}{4} \tr_{g_{\tau}} \left(
\dot{g}_{\tau} \right)
\right|^2 
\leq C,
\label{C2} \\
\left| g_\tau\left(\partial_\tau d \vol^{g_{\tau}}, 
d \vol^{g_{\tau}} \right) \right| 
\leq C,
\label{C3} \\
\begin{aligned}
&\left|
\langle
- \tfrac{1}{2} {\mathfrak D}^{\dot{g}_{\tau}} \phi
+ \tfrac{1}{4} \grad^{g_{\tau}} 
\left( 
\tr_{g_{\tau}} \left( \dot{g}_{\tau} \right)
\right) \bullet_\tau \phi 
- 
\tfrac{1}{4} \div^{g_{\tau}}
\left( \dot{g}_{\tau} \right) \bullet_\tau \phi
, \phi \rangle \right| \\ 
&{\hskip 200pt} \leq C \left(
\left| \n^{g_{\tau}} \phi \right|^2
+ \left| \phi \right|^2
\right), 
\end{aligned}
\label{C4} \\ 
\left| \tfrac{1}{4} \scal^{g_\tau} \right| \leq C.
\label{C5}
\end{gather}
Set $\epsilon = \frac{1}{4C}$ and choose $T>0$ so that 
\begin{equation} \label{assumptiononr}
\frac{1}{4C} \geq \frac{2C^2 + 2C + 3}{4t} + \frac{C}{t^2} 
\end{equation}
for $t > T$.

Take $t > T$. We extend $g_t$ to $t \in \RR$ by setting $g_t =g_0$ 
for $t<0$ and $g_t =g_1$ for $t>1$. Then 
$(\overline{M}_{\infty},\overline{g}_t) = 
(\RR \times M, d\tau^2 + g_{\tau/t})$. 
Take $\phi \in C_0^{\infty} (\Sigma \overline{M}_{\infty})$. 
From (\ref{Diraconcylinder}) we get
\begin{equation*}
\begin{split}
\left| D^{g_{\tau/t}} \phi \right|^2
+
\left| \n_{\nu}^{\overline{g}_t} \phi \right|^2
&=
\left| \left( \nu \cdot D^{\overline{g}_t} - \frac{n}{2} H \right) \phi
\right|^2 \\
&\qquad{} + 
\langle D^{g_{\tau/t}} \phi , \n_{\nu}^{\overline{g}_t} \phi \rangle
+ \langle \n_{\nu}^{\overline{g}_t}\phi, D^{g_{\tau/t}} \phi \rangle.
\end{split}
\end{equation*}
When we integrate over $\overline{M}_{\infty}$ this gives
\begin{equation} \label{firstineq}
\begin{split}
\| D^{g_{\tau/t}} \phi \|^2 
&\leq
2 \| D^{\overline{g}_t} \phi \|^2 + 2 \| \frac{n}{2} H \phi \|^2 \\
&\qquad{} + 
\int_{\overline{M}_{\infty}} 
\left(
\langle D^{g_{\tau/t}} \phi , \n_{\nu}^{\overline{g}_t} \phi \rangle
+ \langle \n_{\nu}^{\overline{g}_t}\phi, D^{g_{\tau/t}} \phi \rangle
\right)
d\vol^{\overline{g}_t}.
\end{split}
\end{equation}
We are going to estimate the terms on the left hand side of this
inequality. 
Define the function $\theta_\tau = g_\tau\left(\partial_\tau
d\vol^{g_{\tau}}, d\vol^{g_{\tau}} \right)$. Then $\partial_\tau
d\vol^{g_{\tau/t}} = \frac{1}{t} \theta_{\tau/t}
d\vol^{g_{\tau/t}}$. For the last term in (\ref{firstineq}) we have
\begin{equation*} 
\begin{split}
&\int_{\overline{M}_{\infty}}
\left(
\langle D^{g_{\tau/t}} \phi , \n_{\nu}^{\overline{g}_t} \phi \rangle
+
\langle \n_{\nu}^{\overline{g}_t}\phi, D^{g_{\tau/t}} \phi \rangle
\right)
d \vol^{\overline{g}_t} \\
&\qquad{}=
\int_{\RR} \int_{ \{ \tau \} \times M} 
\left(
\partial_\tau \langle D^{g_{\tau/t}} \phi , \phi \rangle
- \langle \left[ \n_{\nu}^{\overline{g}_t}, D^{g_{\tau/t}} \right]
\phi , \phi \rangle 
\right)
d \vol^{g_{\tau/t}} d\tau\\
&\qquad{}=
\int_{\RR} \left( \partial_\tau \int_{\{ \tau \} \times M} 
\langle D^{g_{\tau/t}} \phi , \phi \rangle
d \vol^{g_{\tau/t}}
-\int_{\{ \tau \} \times M} 
\langle D^{g_{\tau/t}} \phi , \phi \rangle
\partial_\tau d\vol^{g_{\tau/t}} \right) d\tau \\
&\qquad \qquad{} - \int_{\overline{M}_{\infty}} 
\langle \left[ \n_{\nu}^{\overline{g}_t}, D^{g_{\tau/t}} \right] \phi, 
\phi \rangle d \vol^{\overline{g}_t} \\
&\qquad{}=
- \int_{\overline{M}_{\infty}} 
\left( 
\frac{1}{t} \langle D^{g_{\tau/t}} \phi
, \phi \rangle \theta_{\tau/t} 
+ 
\langle \left[ \n_{\nu}^{\overline{g}_t},
  D^{g_{\tau/t}} \right]\phi , \phi \rangle 
\right)
d \vol^{\overline{g}_t}.
\end{split}
\end{equation*}
so (\ref{firstineq}) becomes
\begin{equation} \label{secondineq}
\begin{split}
\| D^{g_{\tau/t}} \phi \|^2 
&\leq
2 \| D^{\overline{g}_t} \phi \|^2 + 2 \| \frac{n}{2} H \phi \|^2 \\
&\qquad{} - \int_{\overline{M}_{\infty}} 
\left( 
\frac{1}{t} \langle D^{g_{\tau/t}} \phi
, \phi \rangle \theta_{\tau/t} 
+ 
\langle \left[ \n_{\nu}^{\overline{g}_t},
  D^{g_{\tau/t}} \right]\phi , \phi \rangle 
\right)
d \vol^{\overline{g}_t}.
\end{split}
\end{equation}
Since $\partial_\tau (g_{\tau/t}) = \frac{1}{t} \dot{g}_{\tau/t}$
it follows from (\ref{C3}), (\ref{commutator}) and (\ref{C4}) that 
\begin{equation} \label{est2}
\begin{split}
&\left| 
\int_{\overline{M}_{\infty}} \left( \frac{1}{t} 
\langle D^{g_{\tau/t}} \phi, \phi \rangle \theta_{\tau/t} 
+ 
\langle \left[ \n_{\nu}^{\overline{g}_t},
  D^{g_{\tau/t}} \right]\phi , \phi \rangle 
\right) d \vol^{\overline{g}_t}
\right| \\
&\qquad{} \leq
\frac{C}{2t} \left(
\| D^{g_{\tau/t}} \phi \|^2 +  \| \phi \|^2
\right) 
+
\frac{C}{t} \left(
 \| \n^{g_{\tau/t}} \phi \|^2 + \| \phi \|^2
\right).
\end{split}
\end{equation}
By (\ref{C5}) and the Lichnerowicz formula on $(M, g_\tau)$ we have 
\begin{equation} \label{est3}
\begin{split}
\| \n^{g_{\tau/t}} \phi \|^2 
&=
\int_{\RR} \int_{\{ \tau \} \times M} 
\left| \n^{g_{\tau/t}} \phi \right|^2
d\vol^{g_{\tau/t}} dt \\
&=
\int_{\RR} \int_{\{ \tau \} \times M} \left(
\left| D^{g_{\tau/t}} \phi \right|^2
- \tfrac{1}{4} \scal^{g_{\tau/t}} \left| \phi \right|^2
\right)
d\vol^{g_{\tau/t}} dt \\
&\leq
\| D^{g_{\tau/t}} \phi \|^2 + C \| \phi \|^2.
\end{split}
\end{equation}
From (\ref{meancurvature}) and (\ref{C2}) we get 
$\left| \frac{n}{2} H \right|^2 \leq \frac{C}{t^2}$ so
\begin{equation} \label{est1}
\| \frac{n}{2} H \phi \|^2 \leq \frac{C}{t^2} \| \phi \|^2.
\end{equation}
Inserting (\ref{est2}), (\ref{est3}) and (\ref{est1}) into
(\ref{secondineq}) we get
\begin{equation*}
\begin{split}
\| D^{g_{\tau/t}} \phi \|^2
&\leq
2 \| D^{\overline{g}_t} \phi \|^2 + 2 \frac{C}{t^2} \| \phi \|^2 \\
&\qquad{} +
\frac{C}{2t} \left(
\| D^{g_{\tau/t}} \phi \|^2 +  \| \phi \|^2
\right) \\
&\qquad{} + \frac{C}{t}
\left(
\| D^{g_{\tau/t}} \phi \|^2 + C\| \phi \|^2
+ \| \phi \|^2
\right)
\end{split}
\end{equation*}
or 
$$
\| D^{\overline{g}_t} \phi \|^2
\geq
\tfrac{1}{2} \left(1 - \frac{3C}{2t} \right) \| D^{g_{\tau/t}} \phi \|^2
- \left(\frac{C}{t^2} + \frac{C^2 + C}{2t} \right) \| \phi \|^2.
$$
From (\ref{assumptiononr}) we get $1 - \frac{3C}{2t} > 0$ so (\ref{C1}) 
tells us that
\begin{equation*}
\begin{split}
\| D^{\overline{g}_t} \phi \|^2
&\geq
\left(
\frac{1}{2} \left(1 - \frac{3C}{2t} \right) \frac{1}{C}
- \left(\frac{C}{t^2} + \frac{C^2 + C}{2t} \right)
\right)
\| \phi \|^2 \\
&=
\left(
\frac{1}{2C} - \frac{2C^2 + 2C + 3}{4t} - \frac{C}{t^2} 
\right)
\| \phi \|^2 \\
&\geq
\epsilon
\| \phi \|^2.
\end{split}
\end{equation*}
Since $C_0^{\infty} (\Sigma \overline{M}_{\infty})$ is dense in 
$L^2(\Sigma \overline{M}_{\infty})$ we conclude that 
$\overline{g}_t \in \Rinv(\overline{M})$.
\end{proof}
The following corollary is immediate.
\begin{cor} \label{cor:inv=>conc}
Isotopic metrics are concordant.
\end{cor}

\subsection{Surgery}

We are now going to construct a metric with invertible Dirac operator on
the trace of a surgery of codimension $\geq 3$ given such a metric on
the original manifold.

Let $M$ be a closed spin manifold of dimension $n$ and let 
$S^{n-m} \times D^m \to M$ be an embedding. Let $\Sigma$  be the image 
of $S^{n-m} \times \{ 0 \}$. Let $W$ be the trace of the surgery on $M$ 
along $\Sigma$, this can be constructed by attaching 
$D^{n-m+1} \times D^m$ to $M \times [0,1]$ at the image of 
$S^{n-m} \times D^m \times \{ 1 \} \to M \times \{ 1 \}$ and then 
smoothing the corner where the attaching takes place. The trace $W$ is 
a spin manifold with boundary $M \sqcup (\widetilde{M})^-$ where 
$\widetilde{M}$ is the spin manifold obtained from $M$ by surgery
along $\Sigma$. 

\begin{prop} \label{prop_trace}
Assume that $W$ has been constructed from $M$ as above with $m \geq 3$.
Suppose $g \in \Rinv(M)$. Then there is a metric $g^W \in \Rinv(W)$ 
such that $g^W|_{M} = g$.
\end{prop}

The proof is similar to the proof of Theorem 1.2 in 
\cite{Baer_Dahl02}. We need to introduce some notation. 
Suppose $X$ is a submanifold of a Riemannian manifold $Y$.
For $0<r$ define the distance sphere and the distance tube
around $X$ as 
$S_{X}(r) = \{ x \in Y | \operatorname{dist}(x,X) = r \}$
and $U_{X}(r) = \{ x \in Y | \operatorname{dist}(x,X) \leq r \}$.
For $0 < r_1 < r_2$ define the annular region around $X$ as
$A_{X}(r_1,r_2) = \{ x \in Y | 
r_1 \leq \operatorname{dist}(x,X) \leq r_2 \}$.
Let $\nu$ be the outward pointing unit normal of $S_{\Sigma}(r)$ 
and let $dA$ be the volume form of $S_{\Sigma}(r)$.
In \cite[Lemma 2.4]{Baer_Dahl02} the following Lemma is proved in the
case where $X$ is compact, the proof also works in the formulation
here.

\begin{lemma} \label{lemma_BD2.4}
Let $Y$ be a Riemannian spin manifold and let $X \subset Y$ be 
a complete submanifold of codimension $\geq 3$ which has a uniform 
lower bound on the injectivity radius of its normal exponential map 
and for which the second fundamental form of $S_{X}(r)$ is 
bounded for fixed $r$. 

Then there exists $0 < R < 1$ so that for any 
$0 < r < \tfrac{1}{2} R^{11}$ and any smooth spinor field $\phi$ 
defined on $A_{X}(r,(2r)^{1/11})$ satisfying
\begin{itemize}
\item $\int_{S_{X}(\rho)} |\phi|^2 dA$ is finite
for all $\rho \in [r,(2r)^{1/11}]$ and defines a differentiable
function of $\rho$,
\item 
$\int_{S_{X}(\rho)} \operatorname{Re} 
\langle \nabla_{\nu} \phi, \phi \rangle dA$ is finite and non-negative 
for all $\rho \in [r,(2r)^{1/11}]$,
\end{itemize}
it holds that 
$$
\| \phi \|^2_{A_{X}(r,2r)} \leq 
10 r^{5/2}
\| \phi \|^2_{A_{X}( r , (2r)^{1/11} )}.
$$
\end{lemma}

\begin{proof}[Proof of Proposition \ref{prop_trace}]
Since $g \in \Rinv(M)$ there is an $\epsilon^g > 0$ so that 
\begin{equation} \label{surg_lower_bound_orig_metric}
\| D^g \phi \|^2 > \epsilon^g \| \phi \|^2
\end{equation}
for all $\phi \in L^2(\Sigma M)$. 
Proposition 2.1 of \cite{Baer_Dahl02} tells us that there is a
constant $S_0 < 0$ so that for every $S_1 > 0$ there is a metric
$g'$ on $M$ which is conformal to $g$ and has the following
properties: 
\begin{itemize}
\item $g'$ is arbitrarily close to $g$ in the $C^1$-topology on the
space of Riemannian metrics, 
\item $\scal^{g'} \geq S_0$ on all of $M$,
\item $\scal^{g'} \geq 2S_1$ on a neighbourhood $U_0$ of $\Sigma$.
\end{itemize}
The eigenvalues of $D^g$ depend continuosly on the Riemannian metric
with respect to the $C^1$-topology, see for example 
\cite[Prop. 7.1]{Baer96}. We can therefore find a metric $g'$
satisfying the above properties with $S_1 = -8S_0$ while 
(\ref{surg_lower_bound_orig_metric}) holds with the same value of 
$\epsilon^g$. Since $g$ and $g'$ are conformal and the dimension of
the kernel of the Dirac operator is a conformal invariant we get that
$g$ and $g'$ are isotopic and bordant. So if we prove the Theorem for
$g'$ we will also prove it for $g$. We replace our original $g$ with
$g'$.  

Let $r>0$ be a constant so small that 
\begin{itemize}
\item $U_{\Sigma}(2r) \subset U_0$,
\item $(2r)^{1/11} < R$, where $R$ comes from Lemma \ref{lemma_BD2.4}
applied to $\Sigma \subset M$,
\item $(2r)^{1/11} < R$, where $R$ comes from Lemma \ref{lemma_BD2.4}
applied to $\RR \times \Sigma \subset \RR \times M$,
\item $45 r^{1/4} \leq \epsilon^g$.
\end{itemize}
Let $V$ be the trace of the surgery along 
$\Sigma \subset U_{\Sigma}(r)$, this trace is a manifold with boundary
and codimension 2 corners. We divide the boundary of $V$ into a
``horizontal'' part and a ``vertical'' part. The horizontal part
consists of $U_{\Sigma}(r) \sqcup (\widetilde{U})^-$ where
$\widetilde{U}$ is $U_{\Sigma}(r)$ after surgery along  $\Sigma$.
The vertical part is the cylinder 
$[0,1] \times \partial U_{\Sigma}(r)$. The vertical and horizontal
parts meet in the two corners, which are diffeomorphic to 
$\partial U_{\Sigma}(r)$. From \cite{Gajer87} we know that we can 
extend the metric $g$ on $M$ to a metric $g^{V}$ on $V$ 
without decreasing scalar curvature too much. This construction can be
made close to the surgery sphere and we get a metric on $V$ with the
following properties:
\begin{itemize}
\item 
$g^{V}$ is a product metric near the horizontal part of the boundary,
\item 
$g^{V}$ restricts to $g$ on the horizontal part $U_{\Sigma}(r)$ of the
boundary,  
\item
$g^{V}$ restricts to $dx^2 + g$ on a neighbourhood 
$\cong [0,1] \times A_{\Sigma}(r-\delta,r)$ of the vertical part of
the boundary, 
\item $\scal^{g^V} \geq S_1$ on $V$.
\end{itemize}
Define
$(W,g^W) = ([0,1] \times (M - U_{\Sigma}(r)), dx^2 + g) 
\cup (V,g^V)$ where the union is taken along the common boundary
$[0,1] \times \partial U_{\Sigma}(r)$. 

We first prove that 
$\widetilde{g} = g^W|_{\widetilde{M}} \in \Rinv(\widetilde{M})$.
For a contradiction assume that there is a non-trivial harmonic 
spinor field $\phi$ on $(\widetilde{M}, \widetilde{g})$.
Let $\chi: M \to [0,1]$ be a cut-off function with $\chi = 0$ on 
$U_{\Sigma}(r)$, $\chi = 1$ on $M - U_{\Sigma}(2r)$, 
$|\grad \chi| \leq \frac{2}{r}$. Since it has support contained in 
$M - U_{\Sigma}(r) = \widetilde{M} - \widetilde{U}$ 
we can consider $\chi$ also a cut-off function on $\widetilde{M}$.
Set $\psi = \chi \psi$. The spinor field $\psi$ is supported in 
$\widetilde{M} - \widetilde{U}$ where $\widetilde{g} = g$ 
and can be considered a spinor field also for $(M,g)$. 

Since $\scal^{\widetilde{g}} = \scal^{g^V}\geq S_1$ on $\widetilde{U}$ 
and $\scal^{\widetilde{g}} = \scal^{g} \geq 2S_1$ on 
$A_{\Sigma}(r,2r)$ most of the norm of $\phi$ will be concentrated
away from these sets. Lemma 2.2 of \cite{Baer_Dahl02} tells us that 
$$
\| \phi \|^2_{\widetilde{U} \cup A_{\Sigma}(r,2r)} 
\leq 
\frac{-S_0}{S_1 - S_0} \| \phi \|^2_{\widetilde{M}}
=
\frac{1}{9} \| \phi \|^2_{\widetilde{M}}
$$
and it follows that 
\begin{equation} \label{surg_norm_phi_psi}
\begin{split}
\| \psi \|^2_{M} 
&\geq 
\| \psi \|^2_{M - U_{\Sigma}(2r)} \\
&=
\| \phi \|^2_{\widetilde{M} - 
(\widetilde{U} \cup A_{\Sigma}(r,2r))} \\
&\geq 
\tfrac{8}{9} \| \phi \|^2_{\widetilde{M}}.
\end{split}
\end{equation}
Next we are going to show that $\phi$ has even less norm
concentrated in the annular region $A_{\Sigma}(r,2r)$ when compared 
to the larger annular region $A_{\Sigma}(r,(2r)^{1/11})$. 
This will follow from Lemma \ref{lemma_BD2.4} and the fact 
that $\phi$ is harmonic. To apply this Lemma we need to show that
\begin{equation} \label{condition_for_lemma2.4_BD}
\operatorname{Re} \int_{S_{\Sigma}(\rho)} 
\langle \nabla^{\widetilde{g}}_{\nu} \phi, \phi \rangle dA \geq 0
\end{equation}
for all $\rho \in [r,(2r)^{1/11}]$. Choose such a $\rho$ and set 
$\hat{M} = \widetilde{U} \cup A_{\Sigma}(r,\rho)$. Then $\hat{M}$
is a manifold with boundary $\partial \hat{M} = S_{\Sigma}(\rho)$ 
and $\scal^{\widetilde{g}} \geq S_1$ on $\hat{M}$. From the 
Lichnerowicz formula we get
\begin{equation*}
\begin{split}
0
&=
\int_{\hat{M}} \langle (D^{\widetilde{g}})^2 \phi,\phi \rangle 
d\vol^{\widetilde{g}} \\
&=
\int_{\hat{M}} 
\langle (\nabla^{\widetilde{g}})^* \nabla^{\widetilde{g}} \phi,
\phi \rangle d\vol^{\widetilde{g}} 
+
\tfrac{1}{4} \int_{\hat{M}} \scal^{\widetilde{g}} |\phi|^2 
d\vol^{\widetilde{g}} \\
&\geq
\| \nabla^{\widetilde{g}} \phi \|^2_{\hat{M}} 
- 
\int_{\partial \hat{M}} 
\langle \nabla^{\widetilde{g}}_{\nu} \phi, \phi \rangle dA 
+
\tfrac{1}{4} S_1  \| \phi \|^2_{\hat{M}}
\end{split}
\end{equation*}
so 
$$
\operatorname{Re} \int_{\partial \hat{M}} 
\langle \nabla^{\widetilde{g}}_{\nu} \phi, \phi \rangle dA 
=
\int_{\partial \hat{M}} 
\langle \nabla^{\widetilde{g}}_{\nu} \phi, \phi \rangle dA 
\geq
\tfrac{1}{4} S_1 \| \phi \|^2_{\hat{M}}
$$
and (\ref{condition_for_lemma2.4_BD}) follows since $S_1 > 0$.
We now apply Lemma \ref{lemma_BD2.4}, which tells us that 
$$
\| \phi \|^2_{A_{\Sigma}(r,2r)} \leq 
10 r^{5/2}
\| \phi \|^2_{A_{\Sigma}( r , (2r)^{1/11} )}.
$$
Using this estimate we compute
\begin{equation*}
\begin{split}
\| D^g \psi \|^2_M
&=
\| D^{\widetilde{g}} (\chi \phi) \|^2_{\widetilde{M}} \\
&=
\| \grad \chi \cdot \phi \|^2_{\widetilde{M}}\\
&\leq
\frac{4}{r^2} \| \phi \|^2_{A_{\Sigma}(r,2r)} \\
&\leq
40 r^{1/4} \| \phi \|^2_{A_{\Sigma}( r , (2r)^{1/11} )} \\
&\leq
40 r^{1/4} \| \phi \|^2_{\widetilde{M}}
\end{split}
\end{equation*}
which together with (\ref{surg_norm_phi_psi}) and 
the assumption on $r$  tells us that
\begin{equation*}
\begin{split}
\| D^g \psi \|^2_M
&\leq 
45 r^{1/4} \| \psi \|^2_M \\
&\leq 
\epsilon^g \| \psi \|^2_M,
\end{split}
\end{equation*}
and this contradicts (\ref{surg_lower_bound_orig_metric}).

Let $(W_{\infty},g^W)$ be $(W,g^W)$ with half-infinite
cylindrical ends attached. Since $D^{g}$ and $D^{\widetilde{g}}$ 
are both invertible we conclude that the essential spectrum of 
$D^{g^W}$ on $W_{\infty}$ has a gap around $0$, see for example 
\cite[Prop. 3.24]{Bleecker_Booss-Bavnbek03}.
To prove that $g^W \in \Rinv(W)$ it thus remains to show that
$0$ is not an eigenvalue of $D^{g^W}$ on $W_{\infty}$, that is to show
that there are no harmonic spinors in $L^2(\Sigma W_{\infty})$.

To get a contradiction assume that 
$\phi \in L^2(\Sigma W_{\infty})$ is a non-trivial harmonic spinor
field. Then $\phi$ is smooth and the pointwise norm decays
exponentially on the cylindrical ends, see for example 
\cite[Lemma 3.21]{Bleecker_Booss-Bavnbek03}.

Let $V_{\infty}$ be $V$ with the horizontal part of the boundary
extended by half-infinite cylinders. Then 
$(W_{\infty}, g^W) = (\RR \times (M - U_{\Sigma}(r)), dx^2 + g) 
\cup (V_{\infty}, g^V)$ where the union is taken along the common boundary
$\RR \times \partial U_{\Sigma}(r)$. Set 
$\psi = (\chi \circ \pi) \phi$ where $\chi$ is the cut-off
function on $M$ defined above and $\pi: \RR \times M \to M$ is the
natural projection. The spinor field $\psi$ is supported in
$W_{\infty} - V_{\infty} = \RR \times (M - U_{\Sigma}(r))$ where
$g^W = dx^2 + g$ so we can consider $\psi$ to be a spinor field on
$(\RR \times M, dx^2 + g)$.

From \cite[Lemma 2.2]{Baer_Dahl02} applied to 
$V_{\infty} \cup A_{\RR \times \Sigma} (r,2r) \subset  W_{\infty}$
it follows that  
\begin{equation} \label{surg_norm_phi_psi_2}
\| \psi \|^2_{\RR \times M} 
\geq
\tfrac{8}{9} \| \phi \|^2_{W_{\infty}}.
\end{equation}
We now apply Lemma \ref{lemma_BD2.4} to 
$\RR \times \Sigma \subset \RR \times M$. This can be done since 
$|\phi|$ decays exponentially and since the positive scalar curvature 
on $V_{\infty}$ makes the computation for Equation 
(\ref{condition_for_lemma2.4_BD}) work also in this case. The
conclusion is that 
$$
\| \phi \|^2_{A_{\RR \times \Sigma}(r,2r)} \leq 
10 r^{5/2}
\| \phi \|^2_{A_{\RR \times \Sigma}( r , (2r)^{1/11} )}.
$$
Using this, (\ref{surg_norm_phi_psi_2}) and the assumption on $r$ we
compute 
\begin{equation*}
\begin{split}
\| D^{dx^2 + g} \psi \|^2_{\RR \times M}
&=
\| D^{g^W} (\chi \phi) \|^2_{W_{\infty}} \\
&=
\| \grad \chi \cdot \phi \|^2_{W_{\infty}}\\
&\leq
\frac{4}{r^2} \| \phi \|^2_{A_{\RR \times \Sigma}(r,2r)} \\
&\leq
40 r^{1/4} \| \phi \|^2_{A_{\RR \times \Sigma}( r , (2r)^{1/11} )} \\
&\leq
40 r^{1/4} \| \phi \|^2_{W_{\infty}} \\
&\leq 
45 r^{1/4} \| \psi \|^2_{\RR \times M} \\
&\leq
\epsilon^g \| \psi \|^2_{\RR \times M}
\end{split}
\end{equation*}
which is a contradiction since the lower bound
(\ref{surg_lower_bound_orig_metric}) holds also for the product 
Dirac operator $D^{dx^2 + g}$. We conclude 
that the spectrum of $D^{g^W}$ has a gap around $0$, and 
this finishes the proof of the Proposition. 
\end{proof}

\section{Detecting components of $\Rinv(M)$ using the index}
\label{Sec:Detecting}

The alpha invariant of an $n$-dimensional compact spin manifold $M$
without boundary is an element $\alpha(M) \in KO_n(\RR)$ which only
depends on the spin bordism class of $M$. The Index Theorem of Atiyah
and Singer relates the alpha invariant of $M$ to an index-quantity
defined using the kernel of the Dirac operator defined with respect to
some metric. In particular we have the following 
\begin{prop} \label{alpha_is_zero}
Suppose $M$ is a closed spin manifold with a metric $g$ for which $D^g$ is 
invertible. Then $\alpha(M) = 0$.
\end{prop}
The first and obvious conclusion is that $\Rinv(M)$ is empty if 
$\alpha(M) \neq 0$. We are going to use the alpha invariant to
distinguish non-bordant metrics in $\Rinv(M)$ in certain dimensions,
and for this we need some specific manifolds with non-zero alpha 
invariant specified in the following Theorem.
\begin{thm} \label{thm:Carr-Loft-metrics}
For $n = 4k+3$, $k \geq 1$, there are $(n+1)$-dimensional spin
manifolds $Y^i$, $i \in \ZZ$, with boundary $\partial Y^i = S^n$, and
metrics $g^{Y^i} \in \Rpsc(Y^i)$, $i \in \ZZ$, so that 
$\alpha(Y^i \cup_{S^{n}} (Y^j)^-) = c_n (i-j)$ where $c_n \neq 0$.

For $n = 8k$ or $n=8k+1$, $k \geq 1$, there are $(n+1)$-dimensional
spin manifolds $Y^i$, $i =0,1$, with boundary $\partial Y^i = S^n$,
and metrics $g^{Y^i} \in \Rpsc(Y^i)$, $i =0,1$, so that 
$\alpha(Y^1 \cup_{S^{n}} (Y^0)^-) \neq 0$.
\end{thm}

\begin{proof}
In dimensions $n = 4k+3$ manifolds $(Y^i, g^{Y^i})$ with the required
properties are constructed in 
\cite[Ex. 7.6, p. 328]{Lawson_Michelsohn89} using methods of \cite{Carr88}.

For $n = 8k$ and $n=8k+1$ let $Y^0$ be the disc $D^{n+1}$ and let
$g^{Y^0}$ be a positive scalar curvature metric on $Y^0$ which is
equal to the standard metric $g^{S^n}$ on the boundary $S^n$ and is
product in a neighbourhood of the boundary. Let $\Sigma$ be a homotopy
$(n+1)$-sphere with non-vanishing $\alpha$-invariant, see
\cite[Thm. 2.18, p. 93]{Lawson_Michelsohn89}, and let  
$f_0, f_1: D^{n+1} \to \Sigma$ be two disjoint embedded discs. 
Let $W$ be $\Sigma$ with the interiors of $f_0(D^{n+1})$, 
$f_1(D^{n+1})$ removed, then $W$ is a simply connected $h$-cobordism 
with boundary consisting of two components $\partial_0 W = f_0(S^n)$ 
and $\partial_1 W = f_1(S^n)$. By the $h$-Cobordism Theorem there is a
diffeomorphism
$$
(F, \operatorname{id}, f) : 
([0,1] \times \partial_0 W, 
\{ 0 \} \times \partial_0 W, \{ 1 \} \times \partial_0 W)
\to 
(W, \partial_0 W, \partial_1 W).
$$ 
Define $Y^1$ to be $\Sigma$ with the interior of $f_1(D^{n+1})$ 
removed and identify $\partial Y^1$ with $S^n$ using $f_1$. Then  
$Y^1 = W \cup_{\partial_0 W} f_0(D^{n+1})$. On $W$ we set 
$g^{Y^1} = (F^{-1})^* (dx^2 + (f_0^{-1})^* g^{S^n})$
and on $f_0(D^{n+1})$ we set $g^{Y^1} = (f_0^{-1})^* g^{Y^0}$. 
Since $g^{Y^0}$ restricts to $g^{S^n}$ on $S^n$ the definitions 
of $g^{Y^1}$ fit together to a smooth metric of positive scalar
curvature on $Y^1$. Finally 
$$
\alpha(Y^1 \cup_{S^{n}} (Y^0)^-) 
= 
\alpha( ( Z - \operatorname{int} f_1(D^{n+1}) )
\cup_{f_1(S^{n})} f_1(D^{n+1})^- ) 
=
\alpha(Z) \neq 0
$$
and we are done. 
\end{proof}

Define $h^i \in \Rpsc(S^n)$ by $h^i = g^{Y^i}|_{S^n}$.

\begin{thm}
Let $M$ be a compact spin manifold of dimension $n$ and suppose 
$g \in \Rinv(M)$. Then
\begin{itemize}
\item if $n = 4k+3$, $k \geq 1$, there are metrics $g^i \in \Rinv(M)$, 
$i \in \ZZ$, such that $g^i$ is bordant to $g$ and $g^i$ is not 
concordant to $g^j$ for $i \neq j$,
\item if $n = 8k$ or $n=8k+1$, $k \geq 1$,  there is a metric 
$g^1 \in \Rinv(M)$ such that $g^1$ is bordant but not concordant 
to $g$.
\end{itemize}
\end{thm}

\begin{proof}
We prove the Theorem in the case $n = 4k+3$, the other cases are
similar. Fix $i \in \ZZ$. By Proposition \ref{prop_trace} there 
is a metric $g^i$ on $M \# S^n = M$ which is bordant to 
$g \sqcup h^i$ on $M \sqcup S^n$. Since the metric $h^i$ on $S^n$ 
is bordant to the empty manifold by the bordism $(Y^i,g^{Y^i})$ 
we conclude from Corollary \ref{cor:conc&bord_eqrel} that $g^i$ is 
bordant to $g$.

Denote by $(W^i,g^{W^i})$ the bordism between $(M,g^i)$ and
$(M,g)$ we have just constructed. The manifold $W^i$ is 
diffeomorphic to the boundary connected sum of $[0,1] \times M$ 
with $Y^i$.

Take $i,j \in \ZZ$ and suppose the metrics $g^i$ and $g^j$ are 
concordant. By Proposition \ref{prop:attachingboundaries}
we can then find a metric with invertible Dirac operator on
$W^i \cup (W^j)^-$, where the union is obtained by attaching 
the isometric boundary components $(M,g)$ to each other and by 
attaching $(M,g^i)$ to $(M,g^j)$ through a concordance of 
the metrics. Proposition \ref{alpha_is_zero} then tells us that 
$\alpha(W^i \cup (W^j)^-) = 0$. Since $W^i \cup (W^j)^-$ is 
diffeomorphic to the connected sum of $S^1 \times M$ and 
$Y^i \cup_{S^n} (Y^j)^-$ we get 
$0 = \alpha(W^i \cup (W^j)^-) = \alpha(S^1 \times M) 
+ \alpha(Y^i \cup_{S^n} (Y^j)^-) = \alpha(Y^i \cup_{S^n} (Y^j)^-)
= c_n (i-j)$ so $i=j$.
\end{proof}

By Corollary \ref{cor:inv=>conc} this result implies in dimensions 
$n=4k+3$ that if $\Rinv(M)$ is non-empty then it has infinitely many
path-components.  In dimensions $n = 8k, 8k+1$ the result implies that
if non-empty $\Rinv(M)$ has at least two path-components. We conclude
that in these dimensions every closed spin manifold has a metric with
non-trivial kernel of the Dirac operator, which reproves Theorems by
Hitchin \cite[Thm. 4.5]{Hitchin74} and B\"ar \cite[Thm. A]{Baer96}.

\section{Cyclic group actions and metrics with harmonic spinors}
\label{Sec:Cyclic}

Let $M$ be a compact simply connected spin manifold and suppose 
$M \to N$ is a finite covering with covering group $G$. The quotient 
$N$ need not be spin or orientable. We want to find metrics on $M$ 
with non-trivial harmonic spinors. The idea is to use the eta-invariant 
to show that $N$ has metrics with non-trivial harmonic spinors for 
generalized spin structures, and then pull such a metric back to $M$. 
This works under certain conditions on $G$ and $\dim M$, in particular 
we find metrics with harmonic spinors on spheres in all dimensions. 

\begin{thm} \label{exist-main}
Let $M$ be a compact simply connected spin manifold of dimension 
$n \geq 5$ and suppose $M \to N$ is a finite covering with 
covering group $\pi_1(N) = \ZZ/l$. 
\begin{enumerate}
\item If $n$ is odd assume that $N$ is orientable.
\item If $n=2k$ is even assume that $l=2$ and that $N$ is
non-orientable with a $\Pin^+$ structure if $k$ is even or with a
$\Pin^-$ structure if $k$ is odd.
\end{enumerate}
Then $M$ has a $\ZZ/l$-invariant metric with harmonic spinors. 
\end{thm}

The proof relies on work of Botvinnik and Gilkey in 
\cite{Botvinnik_Gilkey96}. Using results of \cite{Gilkey98} and
\cite{Barrera-Yanez99} the argument can also be made to work with
other groups and other assumptions on generalized spin structure on
$N$. The proof will be given through a series of lemmas in the rest of
this section. 

\begin{cor}
For $n \geq 5$ there is a metric with harmonic spinors on the sphere $S^n$.
\end{cor}
\begin{proof}
We obtain a metric with harmonic spinors on $S^n$ by applying Theorem
\ref{exist-main} to the covering $S^n \to P^n$ where $P^n$ is real
projective space of dimension $n$. In odd dimension $P^n$ is
orientable, in even dimension $P^n$ is non-orientable and has a
$\Pin^{\pm}$ structure as required.
\end{proof}

\subsection{Twisted spin structures and Pin structures}
Following \cite{Botvinnik_Gilkey96} we discuss twisted spin 
structures and Pin structures. 

\subsubsection{Twisted spin groups and twisted spin structures}
Let $\ZZ/2$ be the group of two elements written multiplicatively, 
$\ZZ/2 = \{ \pm 1 \}$. Let
\begin{equation} \label{ext_G}
1 \to \ZZ/2 \to {\mathcal G} \overset{\mu}{\to} G \to 1 
\end{equation}
be a central extension of a finite group $G$, this gives an action of
$\ZZ/2$ on ${\mathcal G}$. The group $\Spin(n)$ is a double cover 
$SO(n)$, identifying $\ZZ/2$ with the kernel of the covering homomorphism 
gives an action of $\ZZ/2$ on $\Spin(n)$. Define the twisted spin group 
${\mathcal J}({\mathcal G},\mu,G) = \Spin(n) \times_{\ZZ/2} {\mathcal G}$
where we identify $(\theta,\lambda) = (-\theta,-\lambda)$ for 
$\theta \in \Spin(n)$ and $\lambda \in {\mathcal G}$. The twisted spin 
group ${\mathcal J}({\mathcal G},\mu,G)$ is a double cover of 
$SO(n) \times G$. 

Let $N$ be an $n$-dimensional oriented Riemannian manifold with oriented 
frame bundle $SO(N)$. A ${\mathcal J}({\mathcal G},\mu,G)$ structure on 
$N$ is a principal ${\mathcal J}({\mathcal G},\mu,G)$-bundle 
${\mathcal J}({\mathcal G},\mu,G)(N)$ and an equivariant covering 
${\mathcal J}({\mathcal G},\mu,G)(N) \to SO(N)$ which over open sets $U$ 
in a suitable open cover of $N$ trivializes as 
${\mathcal J} ({\mathcal G},\mu,G) \times U \to SO(n) \times U$. 
A manifold equipped with a ${\mathcal J}({\mathcal G},\mu,G)$ structure 
is called a ${\mathcal J}({\mathcal G},\mu,G)$-manifold. The map $\mu$ 
gives an extension 
$$
1 \to \Spin(n) \to {\mathcal J}({\mathcal G},\mu,G) 
\overset{\mu}{\to} G \to 1, 
$$
through this a ${\mathcal J}({\mathcal G},\mu,G)$ structure 
${\mathcal J}({\mathcal G},\mu,G)(N)$ on $N$ defines 
a homomorphism $\check{\mu}: \pi_1(N) \to G$ as the composition of the 
holonomy of ${\mathcal J}({\mathcal G},\mu,G)(N)$ with $\mu$. 
If $\check{\mu}$ is the trivial 
homomorphism then there is a spin structure $\Spin(N)$ on $N$ so that 
${\mathcal J}({\mathcal G},\mu,G)(N) = \Spin(N) \times_{\ZZ/2} 
{\mathcal G}$.

Suppose $M$ is a compact simply connected spin manifold such that $M$ 
is an oriented covering space of an oriented manifold $N$. 
Let $G = \pi_1(N)$ be the covering group. In 
\cite[Thm. 1.1]{Botvinnik_Gilkey96} a canonical 
${\mathcal J}({\mathcal G},\mu,G)$ structure on $N$ with the property 
that the map $\check{\mu}$ is an isomorphism is constructed. 
The extension $({\mathcal G},\mu,G)$ is given by the lift of the 
action of $G$ on the frame bundle $SO(M)$ to the spin bundle 
$\Spin (M)$ and is split if and only if $N$ is spin. The 
${\mathcal J}({\mathcal G},\mu,G)$ structure on $N$ is given by the 
quotient of $\Spin(M) \times_{\ZZ/2} {\mathcal G}$ by $G$.

\subsubsection{Spinor bundles and Dirac operators for twisted 
spin structures} Let $N$ be a compact oriented $n$-dimensional with a 
${\mathcal J}({\mathcal G},\mu,G)$ structure ${\mathcal J}(N)$. Let $h$ 
be a Riemannian metric on $N$. Let $\alpha$ be a unitary representation 
of ${\mathcal G}$ which is odd with respect to the action of $\ZZ/2$, 
that is $\alpha(-\lambda) = -\alpha(\lambda)$ for all 
$\lambda \in {\mathcal G}$. We denote by 
$\operatorname{Rep}^{\rm odd}({\mathcal G})$ the semi-ring of odd 
unitary representations of ${\mathcal G}$. Let $\Delta$ be the spinor 
representation of $\Spin(n)$, it holds that 
$\Delta(-\theta) = -\Delta(\theta)$ for all $\theta \in \Spin(n)$. Since 
$\Delta(\theta) \otimes \alpha(\lambda) 
= \Delta(-\theta) \otimes \alpha(-\lambda)$ the tensor product 
$\Delta \otimes \alpha$ gives a unitary representation of ${\mathcal J}$. 
Let $\Sigma^{\alpha} N$ be the unitary vector bundle associated to 
${\mathcal J}(N)$ via $\Delta \otimes \alpha$, this is a bundle of 
twisted spinors. As with ordinary spinors there is a Clifford action by 
tangent vectors on $\Sigma^{\alpha} N$, and the Levi-Civita connection 
lifts to a connection on $\Sigma^{\alpha} N$. The Dirac operator 
$D^{h,\alpha}$ acting on sections of $\Sigma^{\alpha} N$ is defined as 
usual.

\subsubsection{Pin groups and Pin structures}
The Clifford algebras $\operatorname{Clif}^{\pm}(n)$ are defined 
as the universal algebra with unit generated by $\RR^n$ with the 
relations $v \cdot w + w \cdot v = \pm 2(v,w)$, $v, w \in \RR^n$.
The groups $\Pin^{\pm}(n)$ are defined as the multiplicative 
subgroups of $\operatorname{Clif}^{\pm}(n)$ generated by the unit
vectors in $\RR^n$. Define $\chi: \Pin^{\pm}(n) \to \ZZ/2$ by
$\chi(v_1 \cdot \dots \cdot v_k) = (-1)^k$ and 
$\Xi^{\pm}: \Pin^{\pm}(n) \to O(n)$ by 
$\Xi^{\pm}(x): v \mapsto \chi(x) x \cdot v \cdot x^{-1}$. Then 
$\Xi^{\pm}$ are two double coverings of $O(n)$ which both restrict
to $\Spin(n) \to SO(n)$.

Let $N$ be an $n$-dimensional Riemannian manifold with frame bundle 
$O(N)$. A $\Pin^{\pm}$ structure on $N$ is a principal 
$\Pin^{\pm}$ bundle $\Pin^{\pm}(N)$ and an equivariant covering 
$\Pin^{\pm}(N) \to O(N)$ which over open sets $U$ in a suitable open 
cover of $N$ trivializes as $\Pin^{\pm} \times U \to O(n) \times U$. 
A manifold equipped with a $\Pin^{\pm}$ structure is called a 
$\Pin^{\pm}$ manifold. If non-empty the set of $\Pin^{\pm}$ structures 
on $N$ is acted on simply and transitively by the cohomology group 
$H^1(N;\ZZ/2)$. So on a simply connected spin manifold there are 
unique $\Pin^{\pm}$ structures given as extensions of the unique 
spin structure.  

\subsubsection{Spinor bundles and Dirac operators for Pin structures}
Let $N$ be a compact $n$-dimensional Riemannian manifold with a 
$\Pin^{\pm}$ structure $\Pin^{\pm}(N)$. Let $\Delta$ be the spinor 
representation of $\Pin^{\pm}(n)$, and let $\Sigma N$ be the unitary 
vector bundle associated to $\Pin^{\pm}(N)$ via $\Delta$, this is
sometimes called a pinor bundle. As with ordinary spinors there is a 
Clifford action by tangent vectors on $\Sigma N$, and the Levi-Civita 
connection lifts to a connection on $\Sigma N$. The Dirac operator $D^h$ 
acting on sections of $\Sigma N$ is defined as usual.

\subsubsection{Pullback to the universal covering space}
Let $N$ be a compact Riemannian manifold with universal covering space
$M$. Assume that $M$ is spin and that $N$ has a ${\mathcal J}$ structure 
${\mathcal J}(N)$ where ${\mathcal J} = {\mathcal J}({\mathcal G},\mu,G)$ 
or ${\mathcal J} = \Pin^{\pm}$. The pullback of ${\mathcal J}(N)$ to 
$M$ is given by an extension of the spin bundle over $M$. In case 
${\mathcal J} = {\mathcal J}({\mathcal G},\mu,G)$ the pullback of 
$\Sigma^{\alpha} N$ is given by $\Sigma M \otimes \CC^d$ where $d$ 
is the dimension of the representation $\alpha$. In case 
${\mathcal J} = \Pin^{\pm}$ the pullback of $\Sigma N$ is given by 
$\Sigma M$. In both cases the pullback of the Dirac operator on $N$ 
defined using some metric is given by the Dirac operator on $M$ with 
the pullback metric.

Let $M$ and $N$ be as in Theorem \ref{exist-main}. If $n$ is odd $N$ has 
a ${\mathcal J}({\mathcal G},\mu,G)$ structure for $G = \ZZ/l$, we say 
that $(N,h)$ has harmonic spinors if $D^{h,\alpha}$ has a non-trivial 
kernel for some 
$\alpha \in \operatorname{Rep}^{\rm odd}({\mathcal G})$.
If $n$ is even we say that $(N,h)$ has harmonic spinors if the Dirac 
operator $D^h$ associated to the $\Pin^{\pm}$ structure has a non-trivial 
kernel. The following Lemma is now obvious.
\begin{lemma} \label{inv_on_quotient}
Let $M$ and $N$ be as in Theorem \ref{exist-main}. If $(N,h)$ has 
harmonic spinors then the pullback of $h$ to $M$ is a 
$\ZZ/l$-invariant metric with harmonic spinors.
\end{lemma}

\subsection{Positive scalar curvature on $N$} 
Using known results on the Gromov-Lawson-Rosenberg conjecture we can 
prove the following Lemma.

\begin{lemma} \label{N_has_PSC}
Let $M$ and $N$ be as in Theorem \ref{exist-main}. 
If $M$ has no $\ZZ/l$-invariant metric with harmonic spinors then $N$ 
has a metric of positive scalar curvature.
\end{lemma}

\begin{proof}
The Gromov-Lawson-Rosenberg conjecture for compact manifolds with
finite fundamental group states the following
\cite[Conj. 5.1]{Rosenberg_Stolz94}: A closed manifold of 
dimension $n \geq 5$ with finite fundamental group admits a metric 
with positive scalar curvature if and only if all index obstructions
associated to Dirac operators with coefficients in flat bundles on
$N$ and its covers vanish. This conjecture is known to be true in the
situation at hand; for orientable manifolds with cyclic fundamental
group by \cite[Thm. 1.1]{Botvinnik_Gilkey97} and
\cite[Thm. A]{Stolz92}, for non-orientable manifolds with fundamental 
group $\ZZ/2$ by \cite[Thm. 5.3]{Rosenberg_Stolz94}. So if $N$ did
not have any metric with positve scalar curvature then the index and
the kernel of some Dirac operator on a cover of $N$ would be non-zero. 
We could then take the pullback of a metric from $N$ to $M$ to produce
a $\ZZ/l$-invariant metric with harmonic spinors on $M$, a contradiction.
\end{proof}

\subsection{The eta invariant}
Let $M$ be a closed Riemannian manifold and let $V$ be a smooth 
vector bundle over $M$. Let $P$ be an operator of Dirac type acting 
on the space of smooth sections of $V$. For complex numbers $z$ with 
large real part the eta function of Atiyah, Patodi and Singer 
\cite{Atiyah_Patodi_Singer75} is defined as 
$\eta(z,P) = \operatorname{Tr}_{L^2}( P (P^2)^{-(z+1)/2})$.
This function has a meromorphic extension to $\CC$ for which $z=0$ is
a regular value and the eta invariant of $P$ is defined as   
$\eta(P) = \half (\eta(0,P) + \dim \ker P)$.

For a closed Riemannian ${\mathcal J}({\mathcal G},\mu,G)$-manifold 
$(N,h)$ and for $\alpha \in \operatorname{Rep}^{\rm odd}({\mathcal G})$
we define $\eta(N,h,\alpha)$ as $\eta (D^{h,\alpha})$.
Let $R^{\rm odd}({\mathcal G})$ be the representation ring associated
to $\operatorname{Rep}^{\rm odd}({\mathcal G})$ and let  
$R^{\rm odd}_0({\mathcal G})$ be the augmentation ideal consisting of
virtual representations of virtual dimension 0. The eta invariant 
$\eta(N,h,\alpha)$ is additive in $\alpha$ so we may extend its 
definition to $\alpha \in R^{\rm odd}({\mathcal G})$.

For a closed $\Pin^{\pm}$ manifold $(N,h)$ we define $\eta(N,h)$ as
$\eta(D^{h})$.

\begin{lemma} \label{eta_ind_of_metric}
Let $M$ and $N$ be as in Theorem \ref{exist-main}.
Let $h^0$, $h^1$ be two metrics on $N$ and assume that $M$ has no
$\ZZ/l$-invariant metric with harmonic spinors.  
\begin{enumerate}
\item If $\dim M$ is odd and $N$ has a 
${\mathcal J}({\mathcal G},\mu,G)$ structure then 
$\eta(N,h^0,\alpha) = \eta(N,h^1,\alpha)$ for all 
$\alpha \in R^{\rm odd}_0({\mathcal G})$.
\item If $\dim M$ is even and $N$ has a $\Pin^{\pm}$ structure then 
$\eta(N,h^0) = \eta(N,h^1)$.
\end{enumerate}
\end{lemma}
\begin{proof}
Let $h_{\tau}$, $\tau \in [0,1]$, be a smooth curve of metrics on $N$
with $h_{\tau} = h^0$ for $\tau$ near $0$ and $h_{\tau} = h^1$ for
$\tau$ near $1$. Lemma \ref{inv_on_quotient} tells us that the Dirac
operator of $h_{\tau}$ is invertible for all $\tau$. Define metrics 
$\overline{h}_t$ on  $\overline{N}_t = [0,t] \times N$ by 
$\overline{h}_t = d\tau^2 + h_{\tau/t}$ for $t > 0$. 
Using the same computation as in Proposition
\ref{prop:stretchcylinder} we conclude that
$(\overline{N}_t, \overline{h}_t)$ has invertible Dirac operator for
$t$ large enough when half-infinite cylinders are attached at the
boundary. 

First suppose that $\dim M$ is odd and that $N$ has a 
${\mathcal J}({\mathcal G},\mu,G)$ structure. Let  
$\alpha \in R^{\rm odd}_0({\mathcal G})$ be the formal
difference of 
$\alpha^+, \alpha^- \in \operatorname{Rep}^{\rm odd}({\mathcal G})$
where $\dim \alpha^+ = \dim \alpha^-$. 
The Atiyah-Patodi-Singer index theorem \cite{Atiyah_Patodi_Singer75}
tells us that 
$$
\ind(D^{\overline{h}_t,\alpha^{\pm}}) 
= 
(\dim \alpha^{\pm}) \int_{\overline{N}_t} \hat{A}(g^{\overline{h}_t})
- 
\epsilon (\eta(N,h^1,\alpha^{\pm}) - \eta(N,h^0,\alpha^{\pm})).
$$
Here $\ind(D^{\overline{h}_t,\alpha^{\pm}})$ is the index of 
$D^{\overline{h}_t,\alpha^{\pm}}$ acting on the space of sections of
the positive half spinor bundle satisfying the Atiyah-Patodi-Singer
boundary condition, $\hat{A}(\overline{h}_t)$ is the $\hat{A}$
differential form computed using the metric $\overline{h}_t$, and
$\epsilon = \pm 1$ is a constant depending only on the dimension. Any
harmonic spinor field satisfying the Atiyah-Patodi-Singer boundary
conditions extends to an $L^2$ harmonic spinor field when
half-infinite cylindrical ends are attached. Since $\overline{h}_t $
has invertible Dirac operator we conclude that the index is zero. We
get 
\begin{equation*}
\begin{split}
\eta(N,h^1,\alpha) - \eta(N,h^0,\alpha)
&=
\eta(N,h^1,\alpha^+) - \eta(N,h^1,\alpha^-) \\
&\qquad{} - \eta(N,h^0,\alpha^+) + \eta(N,h^1,\alpha^-) \\
&=
\epsilon (\dim \alpha^+ - \dim \alpha^-)
\int_{\overline{N}_t} \hat{A}(g^{\overline{h}_t}) \\
&= 0,
\end{split}
\end{equation*}
which proves (1). 

Next suppose that $\dim M$ is even and that $N$ has a $\Pin^{\pm}$
structure. Since $\overline{N}_t$ is then odd-dimensional there is no
integral of a local index density in the index formula for 
$(\overline{N}_t, \overline{h}_t)$, and we have 
$$
\ind(D^{\overline{h}_t}) 
= \epsilon ( \eta(N,h^1) - \eta(N,h^0) ),
$$
where $\epsilon = \pm 1$ is a constant depending only on the
dimension. Again the index vanishes since $\overline{h}_t$ has
invertible Dirac operator and we have proven (2).
\end{proof}

\subsection{Proof of Theorem \ref{exist-main}}
In the work \cite{Botvinnik_Gilkey96} of Botvinnik and Gilkey the
space $\Rpsc(N)$ is studied for a compact manfold $N$ which is either
odd-dimensional with a ${\mathcal J}({\mathcal G},\mu,G)$ structure
and a finite fundamental group satisfying a certain condition or
even-dimensional with fundamental group $\ZZ/2$ and a $\Pin^{\pm}$
structure. The authors construct metrics in $\Rpsc(N)$ with different
values of the eta invariant as follows. Assume $h \in \Rpsc(N)$.
First a (disconnected) manifold $(N',h')$ is found which represents 
zero in an appropriate bordism group and has positive scalar
curvature and non-zero eta invariant. The disjoint union $N \sqcup N'$ 
is then bordant to $N$ and the metric $h \sqcup h'$ of positive scalar 
curvature can be extended over the bordism to give a metric 
$h^1 \in \Rpsc(N)$. The eta invariant is the same for psc-bordant 
metrics so $\eta(N,h^1) = \eta(N,h) + \eta(N,h') \neq \eta(N,h)$.

\begin{proof}[Proof of Theorem \ref{exist-main}]
Assume that $M$ has no $\ZZ/l$-invariant metric with harmonic
spinors. From Lemma \ref{N_has_PSC} we know that $N$ has a metric with
positive scalar curvature. As discussed above the proof of Theorem 3.1 
of \cite{Botvinnik_Gilkey96} gives us two metrics on $N$ with different
$\eta$-invariant, which by Lemma \ref{eta_ind_of_metric} is
impossible.
\end{proof}

\providecommand{\bysame}{\leavevmode\hbox to3em{\hrulefill}\thinspace}
\providecommand{\MR}{\relax\ifhmode\unskip\space\fi MR }
\providecommand{\MRhref}[2]{%
  \href{http://www.ams.org/mathscinet-getitem?mr=#1}{#2}
}
\providecommand{\href}[2]{#2}


\end{document}